\documentclass{article}
\usepackage{amsmath,amssymb}
\usepackage{euscript}
\usepackage{latexsym}
\usepackage[pdftex]{color}
\usepackage[matrix,arrow,curve]{xy}
\usepackage[english]{babel}
\usepackage{wrapfig}
\usepackage{mathrsfs}

\newcommand{\Tor}{{\rm Tor\,}}

\newcommand{\?}{\mathstrut}

\newtheorem{thm}{\bf Theorem}[section]

\newtheorem{corollary}[thm]{\bf Corollary}

\def\ge{\geqslant}

\begin{document}

\begin{center}
{\large
	THE TORSION OF HOMOLOGY GROUPS OF $M(E,I)$-SETS 
}
\\
\medskip
{
Lopatkin V.E
}\\
\end{center}

\section*{Introduction}
In this paper, following \cite{AA,ol}, we consider the homology groups of right pointed sets over a partially commutative monoid $M(E,I)$.\\
\par\textbf{Definition.} Let $E$ be a set and $I \subseteq E \times E$ an irreflexive and summetric relation. A monoid guven by a set of generators $E$ and relations $ab = ba$ for all $(a,b) \in I$ is called \textbf{free partially commutative} \cite{DM} and denoted by $M(E,I)$. If $(a,b) \in I$ then the members $a,b \in E$ are said to be \textbf{commuting generators}.\\
\par The motivation for the research of homology groups of the $M(E,I)$-sets came from a desire to find topology's invariants for asynchronous transition systems. M. Bednarczyk \cite{B} has introduced \emph{asynchronous transition systems} to the modeling the concurrent processes. In \cite{AT} it was proved that the category of asynchronous transition systems admits an inclusion into the category of pointed sets over free partially commutative monoids. Thus asynchronous transition systems may be considered as $M(E,I)$-sets.
\par Follow \cite[Example 3.1]{ol} homology groups of $M(E,I)$-sets may have torsion subgroups. The aim of this paper is to study torsion subgroups of homology groups of followings $M(E,I)$-sets; denote by $X_n^{\bullet}$ a right pointed $M(E,I)$-set, where a set $X_n^{\bullet}$ has $n+2$ elements $X_n^{\bullet}= \{x_0,x_1, \ldots,x_n,*\}$ here $n \ge 1$. Denote by $X_{-1}^{\bullet}=*$ the base point. An action $X_n^{\bullet} \times M(E,I) \to X_n^{\bullet}$ generated by formuls; $x_0 \cdot e = x_1$, $x_1 \cdot e = x_2$,$\dots$, $x_{n-1} \cdot e = x_n$, $x_n \cdot e = *$, $* \cdot e = *$, for all $e \in E$. This $M(E,I)$-set, we'll denote by $X_n^{\bullet}$. \\
\par For any $M(E,I)$-set $X_n^{\bullet}$ let us denote by $(M(E,I)/X_n^{\bullet})^{op}$, or shortly $\mathscr{K}_*(X_m^{\bullet})$ a category which objects are elements of pointed set $X_n^{\bullet}$ and morphisms are triples $(x,\mu,x')$ with $x,x' \in X_n^{\bullet}$ and $\mu \in M(E,I)$ satisfying to $x \cdot \mu = x'$. Sometimes we denote this triples by $x \xrightarrow{\mu} x'$.
\par Let $\mathscr{C}$ be a small category. Denote by $\Delta_{\mathscr{C}}\mathbb{Z}:\mathscr{C} \to \mathrm{Ab}$, or shortly $\Delta \mathbb{Z}$, the functor which has constant values $\Delta \mathbb{Z}(c) = \mathbb{Z}$ at $c \in \mathrm{Ob}\, \mathscr{C}$ and $\Delta \mathbb{Z}(\alpha) = 1_{\mathbb{Z}}$ at $\alpha \in \mathrm{Mor}\,\mathscr{C}$.
\par For any functor $F:\mathscr{C} \to \mathrm{Ab}$, here $\mathscr{C}$ is a small category, denote by $\varinjlim_{k}^{\mathscr{C}} F$ values of the left satellites of the colimit $\varinjlim^{\mathscr{C}}: \mathrm{Ab}^{\mathscr{C}} \to \mathrm{Ab}$. It is well \cite[Prop. 3.3, Aplication 2]{GZ} that thete exists an isomorphim of left satellites of the colimit $\varinjlim^{\mathscr{C}}: \mathrm{Ab}^{\mathscr{C}} \to \mathrm{Ab}$ and the functors $H_n(C_*(\mathscr{C},-)):\mathrm{Ab}^{\mathscr{C}} \to \mathrm{Ab}$. Since the category $\mathrm{Ab}^{\mathscr{C}}$ has enough projectives, these satellites are natural isomorphic to the left derived functor of $\varinjlim^{\mathscr{C}}: \mathrm{Ab}^{\mathscr{C}} \to \mathrm{Ab}$. Denote the values $H_n(C_*(\mathscr{C},-))$ of satellites at $F \in \mathrm{Ab}^{\mathscr{C}}$ by $\varinjlim_n^{\mathscr{C}}F$.
\par We'll be consider any monoid as the small category with the one object. This exert influence on our terminology. In particular a right $M$-set $X$ will be considered and denoted as a functor $X:M^{op} \to \mathrm{Set}$ (the value of $X$ at the unique object will be denoted by $X(M)$ or shortly $X$.) Morphisms of right $M$-sets are natural transformations.

\section{General theorems} Let us show that there is following
\begin{thm}\label{g1}
Let $X_m^{\bullet}$ and $X_{m+1}^{\bullet}$ are $M(E,I)$-sets, here $m \ge -1$, then there exist an isomorphism
$$H_k(X_{m+1}^{\bullet}) \cong H_k(X_m^{\bullet}) \oplus H_{k-1}(E,\mathfrak{M}),$$
here $k \ge 1$ and $H_*(E,\mathfrak)$ is homology groups of simplicial schema $(E,\mathfrak{M})$.
\end{thm}
\textbf{Proof.} Let us consider categories $\mathscr{K}_*(X_m^{\bullet})$ and $\mathscr{K}_*(X_{m+1}^{\bullet})$. By $\mathfrak{In}:\mathscr{K}_*(X^{\bullet}_{m}) \to \mathscr{K}_*(X^{\bullet}_{m+1})$ denote the embedding functor. The functor $\mathfrak{In}$ is defined by formulas $\mathfrak{In}(x_i) = x_{i+1}$ for all $i \in \{0,1,\ldots\,m\}$ on objects, and $\mathfrak{In}(x_j \xrightarrow{e} x_{j+1}) = x_{j+1} \xrightarrow{e} x_{j+2}$, $\mathfrak{In}(* \xrightarrow{e} *) = * \xrightarrow{e} *$ for all $j \in \{0,\dots,m-1\}$ at morphisms. Let us construct a left inverse functor $\mathfrak{Re}:\mathscr{K}_*(X_{m+1}^{\bullet}) \to \mathscr{K}_*(X_m^{\bullet})$ to functor $\mathfrak{In}$. We define the functor $\mathfrak{Re}$ by formulas $\mathfrak{Re}(x_{\xi})=x_{{\xi}-1}$, $\mathfrak{Re}(x_0) = x_0$, $\mathfrak{Re}(*) =*$ for all ${\xi} \in \{1,\dots,m+1\}$ on objects and $\mathfrak{Re}(x_{\xi} \xrightarrow{e} x_{{\xi}+1}) = x_{{\xi}-1} \xrightarrow{e} x_{{\xi}}$, $\mathfrak{Re}(x_0 \xrightarrow{e} x_{1}) = x_0 \xrightarrow{1_{M(E,I)}} x_0$ and $\mathfrak{Re}(* \xrightarrow{e} *) = * \xrightarrow{e} *$ at morphisms. The construction of functors $\mathfrak{In}$ and $\mathfrak{Re}$ is shown in figure \ref{sh}.
\begin{figure}[h!]
$$
  \xymatrix{
 & x_0 \ar@<1ex>@{.>}@/^/[dd]|{\mathfrak{In}} \ar@<1.3ex>[r]^{e_1}  \ar@<-1.3ex>[r]^{\vdots}_{e_s} & x_1 \ar@<1ex>@{.>}@/^/[dd]|{\mathfrak{In}} \ar@<1.3ex>[r]^{e_1}  \ar@<-1.3ex>[r]^{\vdots}_{e_s} & \ldots \ar@<1.3ex>[r]^{e_1}  \ar@<-1.3ex>[r]^{\vdots}_{e_s} & x_m \ar@<1ex>@{.>}@/^/[dd]|{\mathfrak{In}} \ar@<1.3ex>[r]^{e_1}  \ar@<-1.3ex>[r]^{\vdots}_{e_s} & {*} \ar@(r,ur)^(.65){\vdots}_{e_1}  \ar@(dr,r)_{e_s} \ar@<1ex>@{.>}@/^/[dd]|{\mathfrak{In}}\\
&&&&&&\\
x_0  \ar@<1ex>@{.>}@/^/[uur]|{\mathfrak{Re}} \ar@<1.3ex>[r]^{e_1}  \ar@<-1.3ex>[r]^{\vdots}_{e_s} & x_1 \ar@<1ex>@{.>}@/^/ [uu]|{\mathfrak{Re}} \ar@<1.3ex>[r]^{e_1}  \ar@<-1.3ex>[r]^{\vdots}_{e_s} & x_2 \ar@<1ex>@{.>}@/^/ [uu]|{\mathfrak{Re}} \ar@<1.3ex>[r]^{e_1}  \ar@<-1.3ex>[r]^{\vdots}_{e_s} & \ldots \ar@<1.3ex>[r]^{e_1}  \ar@<-1.3ex>[r]^{\vdots}_{e_s} & x_{m+1} \ar@<1ex>@{.>}@/^/ [uu]|{\mathfrak{Re}} \ar@<1.3ex>[r]^{e_1}  \ar@<-1.3ex>[r]^{\vdots}_{e_s} & {*} \ar@(r,ur)^(.65){\vdots}_{e_1}  \ar@(dr,r)_{e_s} \ar@<1ex>@{.>}@/^/ [uu]|{\mathfrak{Re}}
}
$$
\caption{Constructions of functor $\mathfrak{In}:\mathscr{K}_*(X^{\bullet}_{m}) \to \mathscr{K}_*(X^{\bullet}_{m+1})$ and functor $\mathfrak{Re}:\mathscr{K}_*(X^{\bullet}_{m+1}) \to \mathscr{K}_*(X^{\bullet}_{m})$, here the cardinality of $E$ is $s$.}\label{sh}
\end{figure}
\par Denote by $\Delta_{m} \mathbb{Z}$ the functor $\Delta_{\mathscr{K}_*(X^{\bullet}_{m})} \mathbb{Z}:\mathscr{K}_*(X^{\bullet}_{m}) \to \mathrm{Ab}$ for any $m \ge -1$. From \cite[Theorem 3.1]{ol} it follows that homology groups $\varinjlim_n^{\mathscr{K}_*(X^{\bullet}_{m})}\Delta_m \mathbb{Z}$ isomorphic to homology groups of differential object $\left({\?}_m \mathscr{C}_*(\mathscr{K}_*(X^{\bullet}_{m}),\Delta_m \mathbb{Z}),{\?}_m d_* \right)$. Since $\varinjlim_k^{(-)}\Delta \mathbb{Z}$ is exact with respect to first argument, we have homomorphisms; $\mathfrak{In}_k:\varinjlim_k^{\mathscr{K}_*(X^{\bullet}_{m})}\Delta_{m} \mathbb{Z} \to \varinjlim_k^{\mathscr{K}_*(X^{\bullet}_{m+1})}\Delta_{m+1} \mathbb{Z}$ and $\mathfrak{Re}_k:\varinjlim_k^{\mathscr{K}_*(X^{\bullet}_{m+1})}\Delta_{m} \mathbb{Z} \to \varinjlim_k^{\mathscr{K}_*(X^{\bullet}_{m})}\Delta_{m} \mathbb{Z}$, here $k \ge 0$. We have exact sequence of functors
$$ 0 \to \widetilde{ \Delta_{m} \mathbb{Z}} \to \Delta_{m+1} \mathbb{Z} \to {\Delta_{m+1} \mathbb{Z}}/ \widetilde{{\Delta_{m} \mathbb{Z}}}  \to 0,$$
here we've denote by $\widetilde{\Delta_{m} \mathbb{Z}}$ the functor $\Delta_{m} \mathbb{Z}$ with zeros. This exact sequence we can show by following diagram
$$
  \xymatrix{
0 \ar@{->}[d] & 0 \ar@{->}[d] & 0 \ar@{->}[d] & \ldots & 0 \ar@{->}[d] & 0 \ar@{->}[d] \\
0 \ar@<1.3ex>[r] \ar@{->}[d]  \ar@<-1.3ex>[r]^{\vdots} & \mathbb{Z} \ar@{->}[d] \ar@<1.3ex>[r]^{1_{\mathbb{Z}}}  \ar@<-1.3ex>[r]^{\vdots}_{1_{\mathbb{Z}}} & \mathbb{Z} \ar@{->}[d] \ar@<1.3ex>[r]^{1_{\mathbb{Z}}}  \ar@<-1.3ex>[r]^{\vdots}_{1_{\mathbb{Z}}} & \ldots \ar@<1.3ex>[r]^{1_{\mathbb{Z}}}  \ar@<-1.3ex>[r]^{\vdots}_{1_{\mathbb{Z}}} & \mathbb{Z} \ar@{->}[d] \ar@<1.3ex>[r]^{1_{\mathbb{Z}}}  \ar@<-1.3ex>[r]^{\vdots}_{1_{\mathbb{Z}}} & \mathbb{Z} \ar@(r,ur)^(.65){\vdots}_{1_{\mathbb{Z}}}  \ar@(dr,r)_{1_{\mathbb{Z}}} \ar@{->}[d] \\
\mathbb{Z} \ar@{->}[d] \ar@<1.3ex>[r]^{1_{\mathbb{Z}}}  \ar@<-1.3ex>[r]^{\vdots}_{1_{\mathbb{Z}}} & \mathbb{Z} \ar@{->}[d] \ar@<1.3ex>[r]^{1_{\mathbb{Z}}}  \ar@<-1.3ex>[r]^{\vdots}_{1_{\mathbb{Z}}} & \mathbb{Z} \ar@{->}[d] \ar@<1.3ex>[r]^{1_{\mathbb{Z}}}  \ar@<-1.3ex>[r]^{\vdots}_{1_{\mathbb{Z}}} & \ldots \ar@<1.3ex>[r]^{1_{\mathbb{Z}}}  \ar@<-1.3ex>[r]^{\vdots}_{1_{\mathbb{Z}}} & \mathbb{Z} \ar@{->}[d] \ar@<1.3ex>[r]^{1_{\mathbb{Z}}}  \ar@<-1.3ex>[r]^{\vdots}_{1_{\mathbb{Z}}} & \mathbb{Z} \ar@(r,ur)^(.65){\vdots}_{1_{\mathbb{Z}}}  \ar@(dr,r)_{1_{\mathbb{Z}}} \ar@{->}[d] \\
\mathbb{Z} \ar@{->}[d] \ar@<1.3ex>[r]  \ar@<-1.3ex>[r]^{\vdots} & 0 \ar@{->}[d] \ar@<1.3ex>[r]  \ar@<-1.3ex>[r]^{\vdots} & 0 \ar@{->}[d] \ar@<1.3ex>[r]  \ar@<-1.3ex>[r]^{\vdots} & \ldots \ar@<1.3ex>[r]  \ar@<-1.3ex>[r]^{\vdots} & 0 \ar@{->}[d] \ar@<1.3ex>[r]  \ar@<-1.3ex>[r]^{\vdots} & 0  \ar@(r,ur)^(.65){\vdots}  \ar@(dr,r) \ar@{->}[d] \\
0 & 0 & 0 & \ldots & 0 & 0 
}
$$
In this diagram all vertical sequences are exact. Since the functor ${\?}_m \mathscr{C}_*(\mathscr{K}_*(X^{\bullet}_{m}),\Delta_{m} \mathbb{Z})$ is exact with respect to $\Delta_{m} \mathbb{Z}$, we obtain the long exact sequence
$$
  \xymatrix{
\ldots \ar@{->}[r]^(.25){{\partial}_k} & \varinjlim_k^{\mathscr{K}_*(X^{\bullet}_{m})}\Delta_{m} \mathbb{Z} \ar@<1ex>[r]^{\mathfrak{In}_k} & \varinjlim_k^{\mathscr{K}_*(X^{\bullet}_{m+1})}\Delta_{m+1} \mathbb{Z} \ar@<1ex>@{-->}[l]^{\mathfrak{Re}_k}
 \ar@{->}[r]^{\mathfrak{P}_k} & \varinjlim_k^{\mathscr{K}_*(X^{\bullet}_{0})}\mathbb{Z}[x_0] \ar@{->}[r]^(.7){{\partial}_{k-1}} & \ldots 
}
$$
From this long exact sequence we can make up the following shot exact sequence
$$0 \to \mathrm{Im}\,\mathfrak{In}_k \to \underrightarrow{\mathrm{lim}}_k^{\mathscr{K}_*(X_{m+1}^{\bullet})}\Delta_{m+1}\mathbb{Z} \to \mathrm{Im}\,\mathfrak{P}_k \to 0$$
Since $\mathrm{Im}\,\mathfrak{In}_k \cong \varinjlim_k^{\mathscr{K}_*(X^{\bullet}_{m})}\Delta_{m} \mathbb{Z}$, we claim that $\mathrm{Im}\,\mathfrak{P}_k \cong \varinjlim_k^{\mathscr{K}_*(X^{\bullet}_{0})}\mathbb{Z}[x_0]$. Indeed, this follows from long exact sequence and homomorphism's theorem; 
$$\mathrm{Im} \mathfrak{P}_k = \mathrm{Ker}\,\partial_{k-1}, \quad \mathrm{Im}\,\partial_{k-1} = \mathrm{Ker}\,\mathfrak{In}_{k-1} = 0$$
And using homomorphism's theorem, we get $\mathrm{Ker}\,\partial_{k-1} \cong \varinjlim_k^{\mathscr{K}_*(X^{\bullet}_{0})}\mathbb{Z}[x_0]$. Since there exist the homomorphism $\mathfrak{Re}_k:\varinjlim_k^{\mathscr{K}_*(X^{\bullet}_{m+1})}\Delta_{m+1} \mathbb{Z} \to \varinjlim_k^{\mathscr{K}_*(X^{\bullet}_{m})}\Delta_{m} \mathbb{Z}$, we see that shot exact sequence is split, thus we get
\begin{center}
$\varinjlim_k^{\mathscr{K}_*(X^{\bullet}_{m+1})}\Delta_{m+1} \mathbb{Z}
 \cong \varinjlim_k^{\mathscr{K}_*(X^{\bullet}_{m})}\Delta_{m} \mathbb{Z} \oplus \varinjlim_k^{\mathscr{K}_*(X^{\bullet}_{0})}\mathbb{Z}[x_0]
$
\end{center}
using \cite[Theorem~3.1, Example~3.2]{ol}, we'll complete the proof of theorem.\\
\par From this thoerem, we get the following 
\begin{corollary}
Let $X_m^{\bullet}$ and $X_{m+1}^{\bullet}$ are $M(E,I)$-sets, then there exist a isomorphism
$$\Tor (H_k(X_{m+1}^{\bullet})) \cong \Tor (H_k(X_{m}^{\bullet})) \oplus \Tor (H_k(X_{0}^{\bullet})),$$
here $m \ge -1$ and $k \ge 1$.
\end{corollary}
\par For one-dimensional homology groups, we get the following
\begin{thm}
One-dimensional homology groups of $M(E,I)$-sets form $X_m^{\bullet}$ are free, here $m \ge -1$.
\end{thm}
\textbf{Proof.} Indeed, using the theorem \ref{g1}, where $m = 0$, we get the isomorphism $H_1(X_{0}^{\bullet}) \cong H_1(X_{-1}^{\bullet}) \oplus \varinjlim_0^{\mathscr{K}_*(X^{\bullet}_{0})}\mathbb{Z}[x_0]$. From \cite[Example~3.1]{ol}, it follows that the group $H_1(X_{-1}^{\bullet})$ is free. Further, using induction on $m$ and theorem \ref{g1}, we'll complete the proof of this theorem.
\begin{thm}
Let $X^{\bullet} = \{x_0,x_1,\dots,x_n,*\}$ is poset. Suppose that we have a free partially commutative monoid $M(E,I)= M \left(\coprod\limits_{i =1}^{n} E_i, \coprod\limits_{i =1}^{n} I_i \right)$. Let there is an action $X^{\bullet} \times M(E,I) \to X^{\bullet}$ on the pointed set $X^{\bullet}$, thus we get the $M(E,I)$-set $X^{\bullet}$. Suppose that this action generated by formulas, for all $i \in \{1,\dots,n\}$; 
$$x_0 \cdot e^{E_i}_{j_i} = x_1,\,\,x_1 \cdot e^{E_i}_{j_i} = x_2,\dots,x_n \cdot e^{E_i}_{j_i} = *,\,\, * \cdot e^{E_i}_{j_i} = *,$$
here $j_i \in \{j_1,\dots,j_{\mathrm{card} E_i}\}$, and $\mathrm{card}E_i$ is cardinality of $E_i$. Then there exist a isomorphism
$$H_m \left( (M(E,I)/X^{\bullet})^{op} \right)  \cong \bigoplus\limits_{i=1}^{n} H_m \left((M(E_i,I_i)/X^{\bullet})^{op} \right), \qquad m \ge 1.$$
\end{thm}
\textbf{Proof.} Indeed, since $M(E,I)= M \left(\coprod\limits_{i =1}^{n} E_i, \coprod\limits_{i =1}^{n} I_i \right)$, we see that there exist a isomorphism for differentional objects
$$(\mathscr{C}_*((M(E,I)/X^{\bullet})^{op}),d_*) \cong \bigoplus\limits_{i=1}^n ({\?}_i\mathscr{C}_*((M(E_i,I_i)/X^{\bullet})^{op}),{\?}_id_*),$$
but, it knows that the functor $H_*(-,\mathbb{Z})$ is permutable with respect to direct sums. Therefor, we get the isomorphism
$$H_m \left( (M(E,I)/X^{\bullet})^{op} \right)  \cong \bigoplus\limits_{i=1}^{n} H_m \left((M(E_i,I_i)/X^{\bullet})^{op} \right).$$

\end{document}